\documentstyle[12pt]{article}
\newcommand{\qed}{\hfill\rule{4pt}{8pt}\par\vspace{\baselineskip}}
\newtheorem{de}{Definition}[section]
\newtheorem{lm}[de]{Lemma}

\newtheorem{pr}[de]{Proposition}
\newtheorem{co}[de]{Corollary}

\newtheorem{te}[de]{Theorem}

\newtheorem{lemma}[de]{Lemma}

\def\bea{\begin{eqnarray*}}
\def\eea{\end{eqnarray*}}
\def\eps{\varepsilon}

\def\ot{\otimes}
\def\ra{\rightarrow}

\def\phi{\varphi}
\def\ot{\otimes}

\newcommand{\nc}{\newcommand}
\nc{\rnc}{\renewcommand} \nc{\nt}{\newtheorem}

\nc{\thlabel}[1]{\label{theo:#1}}
\nc{\thref}[1]{Theorem~\ref{theo:#1}}
\nc{\selabel}[1]{\label{sect:#1}}
\nc{\seref}[1]{Section~\ref{sect:#1}}
\nc{\lelabel}[1]{\label{lemm:#1}}
\nc{\leref}[1]{Lemma~\ref{lemm:#1}}
\nc{\prlabel}[1]{\label{prop:#1}}
\nc{\prref}[1]{Proposition~\ref{prop:#1}}
\nc{\colabel}[1]{\label{coro:#1}}
\nc{\coref}[1]{Corollary~\ref{coro:#1}}
\nc{\exlabel}[1]{\label{exam:#1}}
\nc{\exref}[1]{Example~\ref{exam:#1}}
\nc{\delabel}[1]{\label{defi:#1}}
\nc{\deref}[1]{Definition~\ref{defi:#1}}
\nc{\remlabel}[1]{\label{rem:#1}}
\nc{\remref}[1]{Remark~\ref{rem:#1}}
\nc{\eqlabel}[1]{\label{equation:#1}}
\nc{\eqref}[1]{(\ref{equation:#1})}

\def\Box{\mbox{$\sqcap\!\!\!\!\sqcup$}}

\begin{document}\title{Hopf Algebras of
Dimension 14}
\author {M. Beattie\thanks{Research supported by NSERC}
\\Mount Allison University\\
Dept. of Mathematics and Computer Science\\Sackville, N.B.Canada E4L 1E6 \\and\\
S. D\u{a}sc\u{a}lescu\thanks{On leave from University of
Bucharest, Facultatea de Matematic\u{a}. Research supported by
Grant SM 10/01 of the Research Administration
of Kuwait University}\\
Kuwait University,
Faculty of Science\\
Dept. Mathematics, PO BOX 5969\\
Safat 13060, Kuwait}
\date{}\maketitle

\begin{abstract}
Let $H$ be a finite dimensional non-semisimple Hopf algebra over
an algebraically closed field $k$ of characteristic 0. If $H$ has
no nontrivial skew-primitive elements, we find some bounds for
the dimension of $H_1$, the second term in the coradical
filtration of $H$. Using these results, we are able to show that
every  Hopf algebra of dimension 14 is semisimple and thus
isomorphic to a group algebra or the dual of a group algebra.
Also   a Hopf algebra of dimension $pq$ where $p$ and $q$ are odd
primes with $p < q \leq  1 + 3p$ and $q \leq 13$ is semisimple
and thus a group algebra or the dual of a group algebra.  We also
have some partial results in the classification problem for
dimension 16.
\end{abstract}

\section{Introduction  }\label{intro}

\paragraph{} In recent years, there has been some progress on the
problem of the classification of finite dimensional Hopf algebras
over an algebraically closed field of characteristic 0. The first
classification discussion appears in \cite{ka} for dimensions 4
and 5, but few techniques were then available. In \cite{williams},
R. Williams classified Hopf algebras of dimension less than 12 by
heavily computational methods.
\paragraph{} Three fundamental results have had a
great impact on the classification theory: the Taft-Wilson Theorem
(\cite{tw}, with a more general version proved in \cite{mont}),
dealing with pointed Hopf algebras, the Nichols-Zoeller Theorem
(\cite{nz}), which extends Lagrange's Theorem to finite
dimensional Hopf algebras, and the Kac-Zhu Theorem (\cite{zhu}),
stating that a Hopf algebra of prime dimension $p$ is isomorphic
to the group algebra $k[C_p]$.   See \cite{and} for a survey of
progress on the classification problem to date.
\paragraph{} It is an open problem to find the smallest dimension
for which infinitely many non-isomorphic Hopf algebras exist. At
 this time, the smallest such dimension where this is known to
happen is 32 (see \cite{grana}, \cite{beattie}). \c{S}tefan showed
that there are only finitely many semisimple Hopf algebras of a
given dimension, up to isomorphism \cite{st}. It has been shown by
a series of results of Masuoka \cite{mas1}, \cite{maspams}, Gelaki
and Westreich \cite{gw} and Etingof and Gelaki \cite{eg} that a
semisimple Hopf algebra of dimension $pq$ where $p$ and $q$ are
  primes is isomorphic to a group algebra or the dual of
a group algebra.
\paragraph{} A new, more conceptual, approach to the
classification of all Hopf algebras in dimensions less than 12
was presented by D. \c{S}tefan in \cite{st2}. The classification
for dimension 12 was done by N. Fukuda in the semisimple case
\cite{fukuda} and completed by S. Natale   \cite{nat}.  In
\cite{an},   N. Andruskiewitsch and S. Natale obtained a series
of general results that led in particular to the classification in
dimensions 15, 21, 25, 35 and 49. Most recently,
 S.-H. Ng \cite{ng} proved that in
dimension $p^2$, $p$ prime, the only Hopf algebras are group
algebras or Taft Hopf algebras. One main aim of this paper is to
prove that in dimension 14 all Hopf algebras are semisimple, so
that the smallest dimension where a Hopf algebra which is neither
semisimple, pointed nor the dual of a pointed, exists is 16.
\paragraph{}  The
classification appears to be  more difficult for even dimensions.
One reason, for instance, is that for odd dimension, if $H$ is a
non-semisimple Hopf algebra, then either $H$ or $H^*$ has a
nontrivial grouplike element \cite{zhu}. The smallest dimensions
that have remained unclassified are 14, 16, 18, 20, etc.

  \paragraph{} We prove the
following.

\begin{te} \label{teorema}
Every Hopf algebra of dimension 14 over an algebraically closed
field $k$ of characteristic 0 is isomorphic to a group algebra or
the dual of a group algebra. Thus there exist three isomorphism
types of $k$-Hopf algebras of dimension 14,   namely $k[C_{14}]$,
$k[D_7]$ and $k[D_7]^*$ where $C_{14}$ is the cyclic group of
order 14 and $D_7$ is the dihedral group of order 14.
\end{te}

To prove  this theorem, we  develop some results about the
dimension of the second term of the coradical filtration of a
finite dimensional Hopf algebra with no nontrivial skew-primitive
elements. As a consequence we are able to prove the following,
from which it follows immediately that a Hopf algebra of dimension
15, 21, 35, 55,  77, 65, 91 or 143 must be isomorphic to a group
algebra or the dual of a group algebra. For dimensions 15, 21,
and 35, this result appears in \cite{an} by different arguments.

\begin{te} \label{teorema2}
  Let $p$ and $q$ be odd primes
where $p < q \leq 1 + 3p$ and $q \leq 13$. Then a  Hopf algebra
$H$ of dimension  $pq$ over an algebraically closed field of
characteristic 0 is semisimple.  Thus, if $H$ has dimension
15,21,35,55,77,65,91 or 143,  $H$ is semisimple and thus
isomorphic to a group algebra or the dual of a group algebra.
\end{te}

 In the
final section, we consider dimension 16, where the classification
of the semisimple Hopf algebras, of the pointed Hopf algebras and
their duals, and of the non-pointed Hopf algebras whose coradical
is a Hopf subalgebra is already done (see \cite{kashina},
\cite{cdr}, \cite{beattie2}, \cite{cdm}).

\section{Preliminaries}
\paragraph{} Except in Section 2, we work over a field $k$ which is
algebraically closed of characteristic 0. All unexplained notation
may be found in \cite{dnr} or in \cite{mont}.  For a coalgebra
$C$, $G(C)$ will denote the group of grouplike elements of $C$.
The coradical filtration of $C$ is $C_0 \subseteq C_1 \ldots$. For
a matrix coalgebra ${\cal M}^c(n,k)$, we will say that $E=\{
e_{ij} |1 \leq i,j \leq n\}$ is a matrix coalgebra basis if $E$ is
a basis for the coalgebra, if $\triangle (e_{ij}) = \sum^n_{r=1}
e_{ir} \otimes e_{rj}$ and if $\epsilon (e_{ij}) = \delta_{ij}$.
\paragraph{}We describe some key notions and results from \cite{st2},
\cite{an} and \cite{nat}.
 First we recall a description of the coradical filtration due to Nichols (unpublished)
 which was key to the results in  \cite{an}. More detail can be found
in \cite[Section 1]{an}.
\paragraph{}For $C$ a coalgebra over an
algebraically closed field $k$, there is a coalgebra projection
$\pi$ from $C$ onto $C_0$ with kernel $I$ \cite[5.4.2]{mont}.
Define $\rho_L: = (\pi \otimes id) \triangle : C \to C_0 \otimes
C$ and $\rho_R: =(id \otimes \pi) \triangle : C \to C \otimes
C_0$. Let $P_n$ be the sequence of subspaces defined recursively
by$$P_0 =0;$$
$$P_1 =\{ x \in C: \triangle (x) =
\rho_L (x) + \rho_R (x) \} = \triangle ^{-1} (C_0 \otimes I + I
\otimes C_0 );$$
$$P_n =\{ x \in C: \triangle (x) - \rho_L (x) - \rho_R (x)
\in \sum_{1 \leq i \leq n-1} P_i \otimes P_{n-i}\}, n \geq 2.$$
Then by a result of Nichols (see \cite[Lemma 1.1]{an}), $P_n =
C_n \cap I$ for $n \geq 0 $. Suppose $C_0 = \oplus_{ \tau \in
{\cal I}} C_{\tau} $ where the $C_{\tau}$ are simple coalgebras
and $C_{\tau}$ has dimension $d^2_{\tau}$. Any $C_0$-bicomodule is
a direct sum of simple $C_0$-subbicomodules and every simple
$C_0$-bicomodule has dimension $d_{\tau}d_{\gamma}$ for some
$\tau, \gamma \in {\cal I}$.

\par Now let $H$ be a finite dimensional Hopf algebra with
coalgebra structure $C$ as above.  Then $P_1$ is a
$C_0$-bicomodule via $\rho_R$ and $\rho_L$.  As in \cite{an} we
denote by  $P_1^{\tau,\gamma}$ the isotypic component of the
$C_0$- bicomodule $P_1$  of type the simple bicomodule with
coalgebra of coefficients $C_{\tau} \otimes C_{\gamma}$.  (If
$C_{\tau} = kg$ for $g$ grouplike, we may use the superscript $g$
instead of $\tau$ in this notation.) Note that $\dim P_1 =
\sum_{\tau,\gamma} \dim P_1^{\tau, \gamma}$ and $\dim H_1 = \dim
H_0 + \dim P_1$. The following lemma is \cite[Corollary 1.3]{an}.
Note that $I, P_1$, etc, are dependent on the choice of the
coalgebra projection $\pi$.

\begin{lm}\label{dimP1} For $S$ the antipode in the Hopf algebra
$H$ and $g \in G(H)$, then $$\dim P_1^{\tau, \gamma} =\dim P_1^{
S\gamma, S\tau} = \dim P_1^{g\tau, g\gamma} = \dim P_1^{\tau g,
\gamma g},$$ where the superscript $S \alpha$ means that the
simple coalgebra is $S(C_{\alpha})$ and the superscript $g \alpha$
or $\alpha g$ means that the coalgebra is $gC_{\alpha}$ or
$C_{\alpha} g$.

\end{lm}

Proofs for the list of results in the next lemma may also be found
in \cite{an}. A nontrivial skew primitive is one which is not
contained in $k[G(H)]$.
 As in \cite{an},  let $H_{0,d}$ with $d \geq1$
denote the direct sum of the simple subcoalgebras of $H$ of
dimension $d^2$.
\begin{lemma}\label{lemma} \begin{enumerate}
\item The order of $G(H)$ divides the dimension of $H_n$ for $n
\geq 0$, and of $H_{0,d}$ for $d \geq 1$.
\item For $P_n$ as
above, $H_n =H_0 \oplus P_n$ and $|G(H)|$ divides $\dim P_n,
\forall n$.
\item Suppose $H$ does not contain any nontrivial
skew primitive element. Suppose that any simple subcoalgebra of
$H$ has dimension 1 or $n^2$, for some $n >1$.Then $n$ divides
$\dim P_1$.
\item Let $H$ be non-semisimple. If $H=H_1$, then $H$ has a non-trivial skew
primitive element.
\item  Let $H$ be non-semisimple. If $\dim H$ is square free, then $H$ has
no non-trivial skew primitive element. \end{enumerate}
\end{lemma}
The next proposition is from \cite{nat} where it is derived from a
result in \cite{st2}.

\begin{pr}\label{natale} Let $H$ be a finite dimensional non-semisimple
Hopf algebra. Suppose that $H$ is generated by a simple
subcoalgebra $C$ of dimension 4 that is stable by the antipode.
Then $H$ fits into an extension
$$1 \to k^G \to H \to X \to 1,$$
where $G$ is a finite group and $X^*$ is a pointed non-semisimple
Hopf algebra.
\end{pr}

Let $H$ be a Hopf algebra of dimension 14. If $H$ is semisimple,
then $H$ is either a group algebra or the dual of a group algebra
by \cite{mas1}. If $H$ is pointed, then by \cite{st1} or by Lemma
\ref{lemma} v), $H$ is semisimple. Thus to prove Theorem
\ref{teorema} we must show that no Hopf algebra of dimension 14
can have coradical containing a matrix coalgebra ${\cal
M}^c(n,k)$ for $n=2$ or 3.

\section{Matrix-like coalgebras} \label{matrixlikecoalgebras}
\paragraph{}In this section $k$ is any field of characteristic different
from 2 but not necessarily algebraically closed. We say that a
coalgebra $C$ is a $2\times 2$ matrix-like coalgebra if $C$ is
spanned by elements $(e_{ij})_{1\leq i,j\leq 2}$ called a
matrix-like spanning set (not necessarily linearly independent)
such that $\triangle (e_{ij})=\sum _{1\leq p\leq 2} e_{ip}\ot
e_{pj}$ and $\epsilon (e_{ij})=\delta _{ij}$ for any $1\leq
i,j\leq 2$. We classify all $2 \times 2$ matrix-like coalgebras
up to isomorphism. This classification can be used in the proof
of Theorem \ref{teorema} but is also of independent interest.
First we describe some $2 \times 2$ matrix-like coalgebras of
dimension less than 4. Let $C_3$ denote  the space with basis $\{
g,h,u\}$ and coalgebra structure defined by
$$\triangle (g)=g\ot
g,\;\; \triangle (h)=h\ot h,\;\; \eps (g)= \eps (h)=1;$$
$$\triangle (u)=g\ot u+u\ot h,\;\; \eps (u)=0.$$
We can take $e_{11}=g, e_{12}=u, e_{21}=0,e_{22}=h$ to be the
matrix-like spanning set.

\paragraph{}For $a \in k$, let $C_2(a)$ be the space with basis $\{ x,y\}$
and coalgebra structure defined by$$\triangle (x)=x\ot x+ay\ot
y,\;\; \eps (x)=1;$$
$$\triangle (y)=x\ot y+y\ot x,\;\; \eps (y)=0.$$
Here we take $e_{11}=x, e_{12}=y, e_{21}=ay,e_{22}=x$.

\paragraph{}A 1-dimensional coalgebra with basis $\{x\}$ is also a $2
\times 2$ matrix-like coalgebra.  Let $e_{11}=x= e_{22},
e_{12}=e_{21}=0$.

\begin{te}\label{matrixlike}
Let $C$ be a $2 \times 2$ matrix-like coalgebra of dimension less
than 4 and with matrix-like spanning set $(e_{ij})_{1\leq i,j\leq
2}$. Then $C$ is isomorphic to one of the following.
\begin{enumerate}\item If $\dim C =3$, then $C \cong C_3$.
\item If $\dim C=2$, then $C \cong C_2 (a)$ for some $a\in k$. If $a\notin
k^2$, then $G(C_2(a))=\emptyset$. If $a\in (k^*)^2$, then
$G(C_2(a))=\{ x+\sqrt{a}y, x-\sqrt{a}y\}$. If $a=0$, then
$G(C_2(a))=\{ x\}$ and in this case $y$ is a $(x,x)$-primitive
element. For $a,b\in k^*$ we have that $C_2(a)\simeq C_2(b)$ if
and only if $\frac{a}{b}$ is a square in $k$. Thus the coalgebras
$C_2(a)$ with $a\neq 0$ are classified by the factor group
$k^*/(k^*)^2$.
\item If $\dim C=1$, then $C$ has basis the group-like element $x$.
 \end{enumerate}\end{te}
{\bf Proof.} We first consider $2\times 2$-matrix-like coalgebras
$C$
of dimension 3.

{\it Case 1.}  The elements $e_{11}, e_{12}, e_{21}$ are linearly
independent. Then write $e_{22}=ae_{11}+be_{12}+ce_{21}$; apply
$\epsilon$ to obtain $a=1$.
 Substitute  $e_{11}+be_{12}+ce_{21}$ for $e_{22}$ in
$\triangle (e_{22})=e_{21}\ot e_{12}+e_{22}\ot e_{22}$ and equate
the coefficients of the basis elements of the tensor product to
see that $1+cb=0$. Thus $b\neq 0$ and $c=-\frac{1}{b}$, which
means that $e_{22}=e_{11}+be_{12}-\frac{1}{b}e_{21}$. Hence we
have that \bea
\triangle (e_{11})&=&e_{11}\ot e_{11}+e_{12}\ot e_{21},\\
\triangle (e_{12})&=&e_{11}\ot e_{12}+e_{12}\ot e_{11}+be_{12}\ot
e_{12}-
\frac{1}{b}e_{12}\ot e_{21},\\
\triangle (e_{21})&=&e_{21}\ot e_{11}+e_{11}\ot e_{21}+be_{12}\ot
e_{21}- \frac{1}{b}e_{21}\ot e_{21}. \eea The dual algebra $C^*$
has dual basis $E_{11},E_{12},E_{21}$ such that $E_{11}$ is the
identity element and$$E_{12}^2=bE_{12}, \quad
E_{21}^2=-\frac{1}{b}E_{21}, \quad
E_{12}E_{21}=E_{11}-\frac{1}{b}E_{12}+bE_{21},\quad
E_{21}E_{12}=0.$$ This shows that there exist two algebra
morphisms from $C^*$ to $k$, namely
$$E_{11}\mapsto 1, E_{12}\mapsto 0, E_{21}\mapsto -\frac{1}{b} \mbox{
and also }E_{11}\mapsto 1, E_{12}\mapsto b, E_{21}\mapsto 0.$$
Therefore $G(C)=\{ e_{11}-\frac{1}{b}e_{21}, e_{11}+be_{12}\}$.
Let $h$ denote $e_{11}-\frac{1}{b}e_{21}$ and let  $g$ denote
$e_{11}+be_{12}$. It is easy to check that  $u=e_{12}$ satisfies
$\triangle (u)=g\ot u+u\ot h$, so then we obtain that $C$ is
isomorphic to $C_3$, as claimed in statement (i).

{\it Case 2.} The elements $e_{11}, e_{12}, e_{22}$ are linearly
independent. Write $e_{21}=ae_{11}+be_{12}+ce_{22}$. By applying
$\eps$, we find $a+c=0$, thus $e_{21}=ae_{11}+be_{12}-ae_{22}$. If
$a\neq 0$, then $e_{11}, e_{12}, e_{21}$ are linearly
independent, and we reduce to Case 1. Assume that $a=0$. Then
$e_{21}=be_{12}$ and if we substitute in $\triangle (e_{21})=
e_{21}\ot e_{11}+e_{22}\ot e_{21}$, we find $be_{11}\ot
e_{12}+be_{12}\ot e_{22}=be_{12}\ot e_{11}+be_{22}\ot e_{12}$,
showing that $b=0$. We obtain \bea \triangle (e_{11})=e_{11}\ot
e_{11}; \quad \triangle (e_{12})=e_{11}\ot e_{12}+e_{12}\ot
e_{22}; \quad \triangle (e_{22})=e_{22}\ot e_{22}, \eea which
again is exactly the coalgebra $C_3$.

\paragraph{}The situation where $e_{11},e_{21},e_{22}$ are linearly
independent reduces to Case 2, and the situation where
$e_{12},e_{21},e_{22}$ are linearly independent reduces to Case 1
by relabelling the spanning elements $e_{ij}$ via the
transposition
$(1\; 2)$. \\[2mm]
We  now consider the situation where $\dim C=2$.

{\it Case 1.} The elements $e_{11}, e_{12}$ are linearly
independent. Then again using the counit, we must have
$e_{21}=ae_{12}$ and $e_{22}=e_{11}+be_{12}$ for some $a,b\in k$.
If $b\neq 0$, then $e_{11},e_{22}$ are linearly independent and
we will investigate this situation in Case 2. Assume $b=0$, thus
$e_{22}=e_{11}$. Since the elements $(e_{ij})_{1\leq i,j\leq 2}$
satisfy the comatrix comultiplication,  \bea \triangle
(e_{11})=e_{11}\ot e_{11}+ae_{12}\ot e_{12} \mbox{ and }
\triangle (e_{12})=e_{11}\ot e_{12}+e_{12}\ot e_{11} \eea
we obtain $C_2(a)$ by letting $x=e_{11}, y=e_{12}$.

\paragraph{}Let $\{ X,Y\}$ be the basis of $C_2(a)^*$ dual to $\{ x,y\}$.
Then $X$ is the identity element of $C_2(a)^*$ and $Y^2=aX$. If
$\phi :C_2(a)^*\ra k$ is an algebra morphism, then $\phi (X)=1$
and $\phi (Y)^2=a$. Hence if $a$ is not a square in $k$, then such
a map $\phi$  cannot exist, and then $G(C_2(a))=\emptyset$. If $a$
is a non-zero square in $k$, then there are two such algebra
morphisms $\phi$, producing two grouplike elements of $C_2(a)$,
namely $x+\sqrt{a}y$ and $x-\sqrt{a}y$. If $a=0$, then
$G(C_2(a))=\{ x\}$, and then clearly $y$ is an $(x,x)$-primitive
element.

\paragraph{}Now let   $a,b\in k^*$ and suppose $C_2(a)\simeq C_2(b)$. This
is equivalent to $C_2(a)^*\simeq C_2(b)^*$, which is equivalent
to the existence of an element $\alpha X+\beta Y\in C_2(a)^*$
such that $\beta \neq 0$ and $(\alpha X+\beta Y)^2=bX$. This is
equivalent to $2\alpha \beta =0$ and $\alpha ^2+\beta^2a=b$.
Since $\beta \neq 0$ we must have $\alpha =0$, and then
$\frac{a}{b}$ is a square in $k$.

{\it Case 2.} The elements $e_{11},e_{22}$ are linearly
independent. Then since $\eps(e_{12})=\eps (e_{21})=0$ we must
have $e_{12}=ae_{11}-ae_{22}, e_{21}=be_{11}-be_{22}$ for some
$a,b\in k$. Then \bea \triangle (e_{11})&=&(1+ab)e_{11}\ot
e_{11}-abe_{11}\ot e_{22}-abe_{22}\ot e_{11}+abe_{22}\ot
e_{22};\\\triangle (e_{22})&=&abe_{11}\ot e_{11}-abe_{11}\ot
e_{22}-abe_{22}\ot e_{11}+(a+1)be_{22}\ot e_{22}; \eea and if we
denote $c=ab, x=\frac{1}{2}(e_{11}+e_{22}), y=e_{11}-e_{22}$,
then $C$ has basis $\{ x,y\}$ and a straightforward computation
shows that
$$\triangle (x)=x\ot x+(c+\frac{1}{4})y\ot y,\;\; \eps (x)=1;$$
$$\triangle (y)=x\ot y+y\ot x,\;\; \eps (y)=0.$$
Thus $C\simeq C_2(c+\frac{1}{4})$.

{\it Case 3.} The elements $e_{11},e_{21}$ are linearly
independent. Then again we see that $e_{12}=ae_{21},
e_{22}=e_{11}+be_{21}$ for some $a,b\in k$. If $b\neq 0$ we
reduce to Case 2. If $b=0$, then $e_{22}=e_{11}$, which means
that $e_{22},e_{21}$ are linearly independent, and we reduce to
Case 2 by relabelling the $e_{ij}$'s via the transposition $(1\;
2)$. All the other cases reduce to one of
these 3 cases via the same relabelling.

Finally, in dimension 1, it is clear  that the coalgebra is a
grouplike coalgebra. \qed

\section{The dimension of $P_1$}\label{P1dim}
\paragraph{}For the remainder of this paper, $k$ will denote an
algebraically closed field of characteristic 0. In this section,
we find lower bounds for the dimension of $P_1$.
   First we prove the following proposition.

\begin{pr}\label{AB} Let $H$ be a finite dimensional
Hopf algebra such that $\dim P_1=n$. If $x_1, x_2, \ldots, x_n$
is a basis for $P_1$, then
$$\triangle (x_i) = \sum^n_{j=1} x_j \otimes A_{ji} + \sum^n _{k=1}
B_{ik} \otimes x_k$$where the $A_{ij}, B_{kl}$ lie in $H_0$ and
satisfy the identities of an $n \times n$ matrix coalgebra. That
is,
$$\triangle (A_{ik}) = \sum ^n_{j=1} A_{ij} \otimes A_{jk} \mbox{ and }
\epsilon (A_{ik}) = \delta_{ik};$$ similarly for the $B_{ik}$. If
$A$ (respectively $B$) is the subcoalgebra of dimension $\leq
n^2$ of $H_0$ generated by the $A_{ij}$ (respectively the
$B_{ij}$) then we cannot have $A
 = B  = k \cdot 1$ or $A \cong B \cong {\cal M}^c (n,k)$.\end{pr}

{\bf Proof.} Since $P_1$ is a right (left) $H_0-$comodule via
$\rho_R$ ($\rho_L$), then the subcoalgebra $A$(respectively $B$)
of $H_0$ satisfies the identities of an $n \times n$ matrix
coalgebra.  If $A = B   = k \cdot 1$, then the $x_i$ are
primitive, contradicting the fact that $H$ is finite dimensional.
If $A \cong B \cong {\cal M}^c (n,k)$, then $P_1$ is a right
$A$-comodule of dimension $n$, and thus is irreducible as a right
$A$-comodule. Thus $P_1$ is irreducible as an $H_0$-bicomodule,
and so $\dim P_1 = n^2$, a contradiction. \qed

\begin{pr}\label{noskewprim}
Let $H$ be a finite dimensional  Hopf algebra with no nontrivial
skew primitive elements. Suppose $H_0 \cong k [G] \oplus {\cal
M}^c (n_1, k) \oplus \ldots \oplus
 {\cal M}^c (n_t, k)$ where $t\geq 1$ and $2 \leq n_1 \leq n_2 \leq \ldots \leq n_t$.
 Then $\dim P_1 $ is a multiple of $|G|$ and is
 greater than or equal to min$(n_1^2, 2n_1|G|)$.\end{pr}

{\bf Proof.} By Lemma \ref{lemma}, $|G|$ divides $\dim P_1$.
Suppose $M$ is a simple $H_0$-sub-bicomodule of $P_1$ so that
$\rho_L(M) \subseteq C_{\tau} \otimes M$ and $\rho_R(M) \subseteq
M \otimes C_{\gamma}$ with $C_{\tau}, C_{\gamma}$ simple
subcoalgebras of $H_0$. If $\dim C_{\tau} = \dim  C_{\gamma} = 1$,
then $H$ contains a nontrivial skew primitive element, a
contradiction. \par If $\dim C_{\tau} = 1$ so that $C_{\tau} = k
\cdot h$ for some $h \in G$, and $\dim C_{\gamma} = n_i^2$ then
$n_i = \dim M \leq  \dim P_1^{h, \gamma}$.  But by Lemma
\ref{dimP1}, $\dim P_1^{h, \gamma} = \dim P_1^{g  h, g \gamma} =
\dim P_1^{S \gamma, h^{-1}} = \dim P_1^{g S \gamma,g h^{-1}}$ for
all $g \in G$, so that, letting $g$ run through $G$, we obtain
$2|G|$ distinct bicomodules $P_1^{\tau,\gamma}$ of dimension
greater than or equal to $n_1$. Thus $\dim P_1 \geq
  2n_1|G|$.\par If $\dim C_{\tau}=n_j$ and $ \dim C_{\gamma}
= n_i$, then dim $P_1 \geq n_in_j \geq n^2_1.$ \qed

\begin{co} \label{cor1} Let $H$ be a finite dimensional  Hopf algebra with
coradical $H_0 \cong k \cdot 1 \oplus {\cal M}^c(n_1, k) \oplus
\ldots \oplus {\cal M}^c(n_t,k)$, where $t \geq 1$, and $2 \leq
n_1 \leq n_2 \leq \ldots \leq n_t$. Then $\dim H > 1+ \sum^t_{i=1}
n^2_i + \min(n_1^2, 2n_1|G|)$.
\end{co}
{\bf Proof.} Since $H$ has no skew primitives, $H \neq H_1$ and so
$\dim H> \dim H_1 = \dim P_1 +1 + \sum^t_{i=1} n_i^2 $. Now apply
Proposition \ref{noskewprim}. \qed

 In
particular, for dimensions 14 and 16, we have the following.
\begin{co} \label{dim} Let $H$ be a non-cosemisimple Hopf algebra.
\begin{enumerate}
\item  If $H$ has dimension 14 or 16, $H$ cannot have coradical
$H_0 = k \cdot 1 \oplus {\cal M}^c(3,k)$. \item If $H$ has
dimension 14 or 16, $H$ cannot have coradical $H_0 = k \cdot 1
\oplus {\cal M}^c(2,k) \oplus {\cal M}^c (2,k)\oplus {\cal M}^c
(2,k)$. \item  If $H$ has dimension 16 or 18, $H$ cannot have
coradical $H_0 \cong k \cdot 1 \oplus {\cal M}^c (2,k) \oplus{\cal
M}^c (3,k)$. \item  If $H$ has  dimension 14, $H$  cannot have
coradical $k[C_2] \oplus {\cal M}^c (2,k) \oplus {\cal M}^c
(2,k)$.
\end{enumerate}
\end{co}
{\bf Proof.}  Statements (i), (ii) and (iii) follow directly from
Corollary \ref{cor1}. Statement (iv) follows from Proposition
\ref{noskewprim} and the fact that since 14 is square-free, by
Lemma \ref{lemma}, the Hopf algebra has no skew primitives. \qed

\section{An approach with injective envelopes}\label{injective}
In this short section, we give a different proof that if $H$ has
dimension 14, then $H$ cannot have coradical $H_0\simeq k[C_2]
\oplus {\cal M}^c(2,k) \oplus {\cal M}^c (2,k)$ and also show that
$H$ cannot have coradical $H_0 \cong k[C_2] \oplus {\cal
M}^c(2,k)$. This approach is also used for other dimensions in
Section \ref{pq}.

\begin{lm} \label{simplefactor} Let $g \in G(H)$, the grouplike
elements of $H$, and let $M$ be a right $H$-subcomodule  of $H$
such that $kg\subset M$ and $M/kg$ is a simple comodule with
coalgebra of coefficients $C$. Then $M\subseteq kg\wedge
C=kg+C+P_1^{g,C}$.
\end{lm}
{\bf Proof.} Let $\{g\} \cup \{x_j|j\in J\}$ be a basis of $M$. If
we replace $x_j$ by $x_j-\eps (x_j)g$, we may assume that $\eps
(x_j)=0$ for any $j$. Let $m\in M$ and write $  \Delta (m)=g\ot
u+\sum _{j\in J}x_j\ot h_j\in M\ot H$. If $\rho :M/kg\rightarrow
M/kg\ot H$ is the comodule structure map of the factor comodule,
and $\hat{m}$ denotes the image of $m \in M$ in the factor
comodule $M/kg$,  then $\rho (\hat{m})=\sum _{j\in J}\hat{x_j}\ot
h_j\in M/kg\ot C$.  Thus $h_j\in C$ for any $j$. Applying $\eps
\ot Id$ to $\Delta(m)$ we see that $u=m$. Thus $\Delta (m)=g\ot
m+\sum _{j\in J}x_j\ot h_j\in kg \otimes M + M \ot C$, and $M
\subseteq kg \wedge C$.  By a result of Nichols, \cite[Lemma
1.2]{an}, $kg\wedge C=kg+C+P_1^{g,C}$. \qed

\begin{co} \label{loewy}
Let $H$ be a non-cosemisimple Hopf algebra with no nontrivial skew
primitives.  Then for any $g \in G(H)$,  there exists a simple
sub-coalgebra $C$ of $H$ of dimension $>1$ such that
$P_1^{g,C}\neq 0$.
\end{co}
{\bf Proof.} Let $E\subseteq H$ be the injective envelope of the
right $H$-comodule $k1$. Since $H$ is not cosemisimple, $k1$ is
not injective (see \cite[Exercise 5.5.9]{dnr}), so $E\neq k1$.
Let $M$ be a subcomodule of $E$ such that $M/k1$ is simple. Then
by Lemma \ref{simplefactor}, $M \subseteq k1 \wedge C =k 1
+C+P_1^{1,C}$ for some simple coalgebra $C$.  Since $H$ has no
nontrivial skew primitives,  $
 P_1^{1,h} = 0$ for $h\in
G(H)$.  Then  $C \neq kh$; otherwise $M\subseteq H_0\cap E=k1$, a
contradiction. Thus $\dim C >1$. Also $P_1^{1,C}\neq 0$,
otherwise the argument above shows that $M\subseteq k1$. Lemma
\ref{dimP1} now implies that $P_1^{g,gC}\neq 0$ for all $g \in
G(H)$. \qed

\begin{co} \label{loewyco} Let $H$ be a finite dimensional
non-cosemisimple Hopf algebra with $H_0 \cong k [G] \oplus {\cal
M}^c(n_1, k) \oplus \ldots \oplus {\cal M}^c (n_t,k)$ with $t$ a
positive integer, $ 2 \leq n_1 \leq n_2 \ldots \leq n_t$, and such
that $H$ has no nontrivial skew primitives. Then
$$\dim H > \dim H_1 = \dim H_0 + \dim P_1 \geq (1 + 2n_1)| G| + \sum^t_{i=1} n^2_i.$$
\end{co}

{\bf Proof.} Let $g\in G$. Then by Corollary \ref{loewy},
$P_1^{g,C}\neq 0$ for some simple sub-coalgebra $C $  of $H$ with
$\dim C  \geq n^2_1$. Then Lemma \ref{dimP1} and the argument in
the proof of Proposition \ref{noskewprim} imply  that $\dim P_1
\geq 2n_1|G|.$ \qed

Corollary \ref{loewyco} may be useful in showing that some Hopf
algebras have skew primitives. For example, let $H$ be a Hopf
algebra of dimension 16 with coradical $ H_0 =k[G] \oplus {\cal
M}^c(2,k)$ with $|G| =4$. Then $H$ has  a nontrivial skew
primitive, and thus contains a pointed Hopf subalgebra of
dimension 8, or else Corollary \ref{loewyco} gives a
contradiction.

The statement in the next corollary  follows easily from Corollary
\ref{loewyco}
  and shows that
the coradical $k[C_2] \oplus \oplus^t_{i=1} {\cal M}^c (2,k)$
cannot occur if $\dim H =14$  for any $t > 0$ , or if $\dim H =
22$ and $t \geq 3$, etc .

\begin{co} \label{casev} Let $H$ be a Hopf algebra of dimension
$pq$, with $p < q$ primes and suppose $H_0 = k[C_p] \oplus
\oplus_{i=1}^t{\cal M}(p,k)$. Then $t < (q - 1 - 2p )/p$. \end{co}

\section{Hopf algebras of dimension 14}\label{main}
In order to prove  Theorem \ref{teorema}, we must show that for a
Hopf algebra $H$ of dimension 14, the  coradical of $H$ cannot be
any of :\begin{enumerate}
\item $H_0=k\cdot 1\oplus C$, with $C\simeq {\cal M}^c(2,k)$ .
\item $H_0=k\cdot 1\oplus C$, with $C\simeq {\cal M}^c(3,k)$.
\item $H_0=k\cdot 1\oplus C\oplus D$, with $C\simeq D\simeq {\cal M}^c(2,k)$.
\item $H_0=k\cdot 1\oplus C\oplus D\oplus E$, with
$C\simeq D\simeq E\simeq {\cal M}^c(2,k)$.
\item $H_0=k[C_2]\oplus C$, with $C\simeq {\cal M}^c(2,k)$.
\item $H_0=k[C_2]\oplus C\oplus D$, with $C\simeq D\simeq {\cal M}^c(2,k)$.
\item $H_0=k[C_2]\oplus C$, with $C\simeq {\cal M}^c(3,k)$.
\item $H_0=k[C_7]\oplus C$, with $C\simeq {\cal M}^c(2,k)$.\end{enumerate}
By Lemma \ref{lemma}, cases (vii) and (viii) cannot occur.
 Cases (ii), (iv) and (vi) were eliminated by Corollary \ref{dim}
 and Case (v) by Corollary \ref{casev}.
  It still remains to consider cases (i) and (iii).
  The key to the argument here is Proposition \ref{natale}, which
yields the next lemma directly.
\begin{lm}  \label{lema3.6}
No non-semisimple Hopf algebra $H$ of squarefree dimension can be
generated as a Hopf algebra by a simple subcoalgebra of dimension
4 that is stable under the antipode.\end{lm}
{\bf Proof.} If $H$
were generated as a Hopf algebra by a simple subcoalgebra of
dimension 4 that is stable under the antipode, then by
Proposition \ref{natale}, $H$ fits into an extension
$$1 \to k^G \to H \to X \to 1,$$
where $G$ is a finite group and $X^*$ is a pointed non-semisimple
Hopf algebra. Then $X^*$ contains non-trivial skew primitive
elements, and then so does $H^*$. Now we obtain a contradiction by
Lemma \ref{lemma}(v). \qed

\begin{co} \label{casei}
For $p$ and $q$ distinct primes, no  Hopf algebra $H$ of
dimension $pq$ can have coradical $H_0$ isomorphic to $k[G] \oplus
C$  where $C \cong {\cal M}^c(2,k)$ and $G = G(H)$.
\end{co}
{\bf Proof.}  Since the antipode maps simple subcoalgebras to
simple subcoalgebras, $S(C)=C$. Then the Hopf subalgebra $L$ of
$H$ generated by $C$ must be all of $H$. For, by the
Nichols-Zoeller Theorem, $\dim L$ is $p, q$ or $pq$ and by the
Kac-Zhu Theorem, the only Hopf algebras of prime dimension are the
group algebras.
 \qed

Thus Case (i) (and also Case (v)) cannot occur. To eliminate Case
(iii), it suffices to show that $H$ contains a simple
subcoalgebra $C$ of dimension 4 stable under $S$. Then $C$
generates a  Hopf subalgebra which again must be all of $H$.

\begin{pr} \label{stable} Let $H$ be a  Hopf algebra  such that $S^4 =Id$ and
$H_0=k \cdot 1\oplus C\oplus D$ with $C\simeq D\simeq {\cal
M}^c(2,k)$.  Then $H$ contains a 4-dimensional simple coalgebra
stable under the antipode.\end{pr}
{\bf Proof.} If $S(C)=C$, then
  we are done, so we need to deal with
the case where $S(C) =D$ and $S(D) =C$. By \cite[Theorem 1.4
(b)]{st2}, we may choose a  matrix coalgebra basis $\{e_{ij} | 1
\leq i,j\leq 2\}$ for $C$ such that
$$S^2 (e_{ij}) =  (-1)^{i+j} e_{ij}.$$
Denote $S(e_{ij})$ by $f_{ji}$. Then $\{ f_{ij}| 1 \leq i,j \leq
2\}$ is a matrix coalgebra basis for $D$. Note that $S(f_{ji}) =
(-1)^{i+j} e_{ij}.$ Since $\epsilon (e_{12}) = \epsilon (e_{21})
=0$, we have that $m(S \otimes Id) \triangle (e_{ij}) = m (Id
\otimes S) \triangle (e_{ij}) =0$ for $i \neq j$ and thus:
\begin{equation}\label{eqn1}f_{11} e_{12} + f_{21} e_{22} =
 e_{11} f_{21} + e_{12} f_{22} = e_{21} f_{11} + e_{22} f_{12} =
 f_{12} e_{11} + f_{22} e_{21} =0. \end{equation}
Similarly, since $\epsilon (f_{ij}) =0$ for $i \neq j$, we have
$m(S \otimes Id) \triangle (f_{ij}) = m(Id \otimes S) \triangle
(f_{ij}) =0$ for $i \neq j$ and thus:
\begin{equation}\label{eqn2} e_{11} f_{12} - e_{21} f_{22} =
 f_{12} e_{22} - f_{11} e_{21} = e_{22} f_{21} - e_{12} f_{11} =
 f_{21} e_{11} - f_{22} e_{12} =0. \end{equation}
Now let $E_{11} = e_{11} f_{22}, E_{12} = e_{12}f_{21}, E_{21} =
e_{21} f_{12}, E_{22} = e_{22} f_{11}$. Then $\epsilon (E_{ij}) =
\epsilon (F_{ij}) = \delta_{ij}$ for $1 \leq i,j \leq 2$. Also
$$\triangle E_{11} =
E_{11} \otimes E_{11} + E_{12} \otimes E_{21} + e_{11} f_{21}
\otimes e_{11} f_{12} + e_{12} f_{22} \otimes e_{21} f_{22}$$ and
by equations (\ref{eqn1}) and (\ref{eqn2}), the sum of the last
two terms is 0. Similar computations show that $\triangle E_{ij}
= \sum^2_{k=1} E_{ik} \otimes E_{kj}$, and since $S(E_{ii}) =
E_{jj}$ and $S(E_{ij}) = - E_{ij}$ for $i \neq j$, the coalgebra
$E$ generated by the $E_{ij}$'s is invariant under $S$. If this
coalgebra is a  4-dimensional coalgebra, then it is simple and we
are done. Similarly let $F_{ii} =f_{ii} e_{jj}, F_{ij} = f_{ij}
e_{ji}, i \neq j$. Then $\triangle F_{ij} = \sum^2_{k=1} F_{ik}
\otimes F_{kj}$ and let $F$ be the coalgebra generated by the
$F_{ij}$. Again, if the $F_{ij}$ are linearly independent, then
$F$ is a 4-dimensional simple coalgebra, invariant under $S$, and
we are done. Suppose $\dim E <4$ and $\dim F<4$. If $\dim E$ or
$\dim F$ is 2 or 3, then by the results of Section 2, $H$ contains
a nontrivial grouplike element, which is impossible.  Thus $E = F
= k \cdot 1$. Then $E_{ij} = F_{ij} =0$ for $i \neq j$ and $E_{ii}
= F_{ii} =1$, i.e.
\begin{equation}\label{eqn3} e_{11} f_{22} = f_{22} e_{11} =
e_{22} f_{11} = f_{11} e_{22} = 1. \end{equation}
\begin{equation}\label{eqn4} e_{12} f_{21} = f_{21} e_{12} =
e_{21} f_{12} = f_{12} e_{21} =0. \end{equation} Now $m(S \otimes
Id) \triangle (e_{ii}) =1 = m(Id \otimes S) \triangle (e_{ii})$
yields
\begin{equation}\label{eqn5} 1= e_{22} f_{22} + e_{21} f_{21} \end{equation}
\begin{equation}\label{eqn6} 1= e_{11} f_{11} + e_{12} f_{12} =
f_{11} e_{11} + f_{21} e_{21}. \end{equation} Thus \bea
f_{12} = (f_{11} e_{11} + f_{21} e_{21}) f_{12} \mbox{ by (\ref{eqn6}) }
= f_{11} e_{11} f_{12} \mbox{ by (\ref{eqn4})} \eea and so $x
= e_{22} f_{12} = e_{11} f_{12} = e_{21} f_{22} = -e_{21} f_{11}$
by (\ref{eqn3}), (\ref{eqn2}) and (\ref{eqn1}).
Now \begin{eqnarray*} \triangle (x) &=& \triangle (e_{11}f_{12})\\
&=& e_{11} f_{11} \otimes e_{11} f_{12} + e_{11} f_{12} \otimes
e_{11} f_{22} + e_{12} f_{11} \otimes e_{21} f_{12} + e_{12}
f_{12} \otimes e_{21}f_{22} \\
&=& (e_{11} f_{11} + e_{12}f_{12}) \otimes x + x \otimes
e_{11}f_{22} \mbox{ by (\ref{eqn3}) and (\ref{eqn4}) }\\
&=& 1 \otimes x + x \otimes 1  \mbox{ by (\ref{eqn6}) }.
\end{eqnarray*}
Therefore $0= x = e_{22} f_{12}$ and since $e_{22}$ is a unit,
$f_{12}=0$ which contradicts the fact that the $f_{ij}$ are a
basis for a 4-dimensional subspace.  \qed
\begin{co} A  Hopf algebra of dimension 14 cannot have coradical
$k \cdot 1 \oplus {\cal M}^c(2,k) \oplus {\cal M}^c (2,k).$
\end{co} {\bf Proof.}  We can assume that $H^*$ does not contain
any nontrivial grouplike. Indeed, if $H^*$ had a nontrivial
grouplike, then $H^*$ would be in one of the situations
(v)-(viii), which  have been shown to be impossible. Therefore
since by \cite{radford} the order of $S$ divides $4(lcm (|G(H)|,
|G(H^*)|)$ $ =4,$ and since $H$ is not semisimple, we have $S^2
\neq Id$ and $S^4 = Id$ and the statement follows from Proposition
\ref{stable}.\qed This completes the proof of Theorem
\ref{teorema}.

\section{Application to odd dimension $pq$}\label{pq}

\paragraph{} Now
we apply the previous results to prove Theorem \ref{teorema2},
i.e., we  show that a Hopf algebra $H$ of dimension $pq$ where $p
$ and $q$ are odd primes with $p < q \leq 1 + 3p$ and $q \leq 13$
is semisimple. \par This gives another proof that Hopf algebras of
dimension 15, 21 and 35 are semisimple (see \cite{an}) and adds
the fact that Hopf algebras of dimension $55 , 77 , 65 , 91 \mbox{
and } 143 $ are semisimple. The proof depends upon the techniques
developed in the previous sections and on the following lemmas.

\begin{lm} \label{gplike} \cite {zhu}
If $H$ is of odd dimension and is not semisimple,
then $H$ or $H^*$ has a nontrivial grouplike. \end{lm}

\begin{lm} \label{Cp} Let $H$ be a nonsemisimple Hopf algebra
with a nontrivial grouplike element and dimension $pq$ where $p$
and $q$ are odd primes with $p <q$. Then $G(H) \cong C_p$ and the
antipode $S$ of $H$ has order $4p$.
\end{lm}

\noindent{\bf Proof.}  The statement follows directly from \cite
[Proposition 5.2]{ng}.   \qed

Now we prove our claim about Hopf algebras of dimension $pq$.\\

\noindent {\bf Proof of Theorem \ref{teorema2}.} Suppose $H$ is
not semisimple and thus has no nontrivial skew primitives by Lemma
\ref{lemma} (v). By Lemma \ref{gplike}, without loss of
generality, we may assume that $H$ has a nontrivial grouplike.
Then by Lemma \ref{Cp}, $G(H) \cong C_p$. Suppose $H_0 \cong
k[C_p] \oplus {\cal M}^c(n_1,k) \oplus \ldots \oplus {\cal M}^c
(n_t,k)$ where $n_1 \leq n_2 \leq \ldots \leq n_t$.
\par Suppose
 $n_i \geq p$ for all $i$.   Then by Corollary \ref{loewyco}, we have $\dim H = pq
> (1 + 2p)p + p^2=p+3p^2$ so that $q > 3p + 1$ which is a
contradiction.\par
 Now suppose $n_1<p$. Here several cases occur.
\par
Case i) Suppose $H_0 \cong k[C_p] \oplus \oplus^p_{i=1}D_i$ where
$D_i \cong {\cal M}^c(2,k)$. Let ${\cal D} =\{ D_1, \ldots,
D_p\}.$ By the argument in the proof of Corollary \ref{casei}, we
cannot have $S(D_i) =D_i$ for any $i$. Thus $S$ induces a
permutation $\sigma$ on the set ${\cal D}$ with no fixed points.
The order of $\sigma$ is a divisor of $4p$, so any cycle has
length 2,4, or $p$ in the decomposition of $\sigma$ as a product
of disjoint cycles.  If there were a cycle of length 2 or 4, then
$\sigma$ has only cycles of even length which is impossible since
$p$ is odd. Thus $\sigma$ is a cycle of length
$p$ and
we may assume that $S(D_i) = D_{i+1}$, with subscripts modulo $p$.
\par By the argument in the proof of Corollary \ref{loewyco}, and
using the obvious notation for the simple subcoalgebra, we have
that for all $g\in G \cong C_p$, $\dim P_1^{g,i} \geq 2$ for some
$i$. Then also,   since by Lemma \ref{dimP1},
$$\dim P_1^{Si, g^{-1}}= \dim P_1^{g, S^2i} = \ldots = \dim
P_1^{i,g^{-1}} = \dim P_1^{g,Si} = \dim P_1^{S^2i, g^{-1}} =
\ldots = \dim P_1^{g,i}  $$ then $\dim(\oplus^{2p-1}_{i=0} S^i
(P_1^{g,i})) \geq 4p$. Since this inequality holds for every $g
\in G$, we have that $\dim P_1 \geq 4p^2$. But then $pq > 5p
+4p^2$ so that $q >5 +4p$. But then $5+4p <q\leq 1+3p$, which is a
contradiction.
\par
Case ii) Suppose $H_0 \cong k[C_p] \oplus \oplus^p_{i=1}D_i \oplus
\oplus^t_{j=1} {\cal M}^c (n_j, k)$ where $D_i \cong {\cal
M}^c(2,k)$, $n_j \geq 2$ and $t \geq 1$.  If $n_1 < p$, then $\dim
H_0 \geq 5p + pn_1^2 \geq 9p$ and since $\dim P_1 \geq 4p$, the
inequality $pq
> 13p$ together with the condition that $q \leq 13$ gives a
contradiction.  If $n_1 \geq p$ then $\dim H_0 \geq 5p + p^2$,
$\dim P_1 \geq 4p$ and thus $q > 9 + p$.  Since $p \geq 3$, this
means that $q > 12$ and since $q$ is prime, $q \geq 13$.  Now if
$q = 13$, then $p \geq 5$ and we do not have $13 > 9 + 5 = 14$.
Thus this case is impossible.
\par
Case iii) Suppose $H_0 \cong k [C_p] \oplus \oplus^p_{i=1} {\cal
M}^c (n,k) \oplus^t_{j=1} {\cal M}^c (n_j, k)$ with $2< n <p , t
\geq 0.$  Then $\dim H_0 \geq 10p$ and $\dim P_1 \geq 6p$ so $pq
> 16 p$ and $q > 16$, a contradiction.

\begin{co} Hopf algebras of dimensions 143, 91, 65,  55, 77,
21, 35, 15 are semisimple.
\end{co}

\section{Hopf algebras of dimension 16}
\paragraph{} The next
dimension  for which the classification is incomplete is
dimension 16. All pointed Hopf algebras of dimension 16 were
described in \cite{cdr} and semisimple Hopf algebras of dimension
16 were classified in \cite{kashina}. All non-pointed
non-semisimple Hopf algebras of dimension 16 with coradical a
Hopf subalgebra are classified in \cite{cdm}.  This leaves the
problem of classifying Hopf algebras of dimension 16 whose
coradical is not a Hopf subalgebra.

\begin{pr} A Hopf algebra of dimension 16 cannot have coradical
isomorphic to any of\begin{enumerate}\item $k \cdot 1 \oplus
{\cal M}^c(2,k)$
\item $k \cdot 1 \oplus {\cal M}^c(2,k)\oplus {\cal M}^c(2,k)\oplus
{\cal M}^c(2,k)$\item $k \cdot 1 \oplus {\cal M}^c(3,k)$
\item $k \cdot 1 \oplus {\cal M}^c(2,k)\oplus {\cal M}^c(3,k)$\end{enumerate}
\end{pr}{\bf Proof.} Impossibility of the last three coradicals follows from
Corollary \ref{dim}. Suppose that $H$ is a Hopf algebra of
dimension 16 with coradical $H_0 = k \cdot 1 \oplus C$ with $C
\cong {\cal M}^c(2,k)$. Then $C$ is $S$-stable and generates a
Hopf subalgebra $L$. If $L=H$, then by Proposition \ref{natale},
$H$ fits into an extension
$$1 \to k^G \to H \to X \to 1.$$
Then the dimension of $X$ is greater than or equal to 4, since
$X^*$ is pointed and nonsemisimple, so the order of $G$ is 1, 2 or
4. If $|G| =1$, then $H^*$ is pointed and this is impossible by
the classification of the duals of the pointed Hopf algebras of
dimension 16 in \cite{beattie2}. If $|G| =2$ or 4 then $k^G$ is
isomorphic to a group algebra and $H$ has nontrivial grouplikes.
Suppose $L$ has dimension 8. Since $H_0$ has dimension 5, $L$ is
not cosemisimple and so $L$ is the unique nonpointed
noncosemisimple Hopf algebra of dimension 8 (see \cite{st2}). But
then $G(L) \cong C_2$, so this is also impossible. \qed

It remains to eliminate $k \cdot 1 \oplus {\cal M}^c(2,k) \oplus
{\cal M}^c(2,k)$ to conclude that nonsemisimple Hopf algebras of
dimension 16 have a nontrivial grouplike. By Lemma \ref{lemma},
coradicals containing $k[G] \oplus {\cal M}^c(3,k), |G| = 2$ or 4,
cannot occur. So the coradical cannot contain a copy of ${\cal
M}^c(3,k)$. It is shown in \cite{beattie2} that $k[C_2] \oplus
\oplus^t_{i=1} {\cal M}^c(2,k)$ with $t =1,2,3$ all occur as
coradicals of duals of pointed Hopf algebras of dimension 16, but
it is not known whether every Hopf algebra with such a coradical
has a pointed dual.

\vspace{1.0in}

\begin{center} {\bf Acknowledgement} \end{center}
\par The second author would like to thank Mount Allison University
for their kind hospitality during his visit there in June, 2001.
\par The authors thank the referee for suggesting shorter proofs
for Propositions \ref{AB} and \ref{noskewprim}, which then led to
a sharpening of the inequalities in Proposition \ref{noskewprim}
and Corollary \ref{loewyco} and to an improvement of Theorem 0.2.

\end{document}